\documentclass[a4paper,11pt]{article}
\usepackage{amsfonts}
\usepackage{amstext}
\usepackage{amsthm}
\usepackage{amsmath}
\usepackage{amssymb}
\usepackage{latexsym}
\usepackage[latin1]{inputenc}
\newcommand{\bl}{\hfill\rule{2mm}{2mm}}
\newcommand{\R}{\Bbb{R}}

\newtheorem{teor}{Theorem}[section]

\newtheorem{cor}{Corollary}[section]

\newcommand{\n}{\noindent}

\begin{document}

\title{On the geometric dependence of Riemannian Sobolev best constants
}
\author{\textbf{Ezequiel R. Barbosa, Marcos Montenegro
 \footnote{\textit{E-mail addresses}:
ezequiel@mat.ufmg.br (E. R. Barbosa), montene@mat.ufmg.br (M.
 Montenegro)}}\\ {\small\it Departamento de Matem\'{a}tica,
Universidade Federal de Minas Gerais,}\\ {\small\it Caixa Postal 702,
 30123-970, Belo Horizonte, MG, Brazil}}
\date{} \maketitle

\markboth{abstract}{abstract}
\addcontentsline{toc}{chapter}{abstract}

\hrule \vspace{0,2cm}

\n {\bf Abstract}

We concerns here with the continuity on the geometry of the second Riemannian $L^p$-Sobolev best constant $B_0(p,g)$ associated to the AB program. Precisely, for $1 \leq p \leq 2$, we prove that $B_0(p,g)$ depends continuously on $g$ in the $C^2$-topology. Moreover, this topology is sharp for $p = 2$. From this discussion, we deduce some existence and $C^0$-compactness results on extremal functions.

\vspace{0.5cm}
\hrule\vspace{0.2cm}

\section{Introduction and main results}

Best constants and sharp first-order Sobolev inequalities on compact
Riemannian manifolds have been extensively studied in the last few
decades and surprising results have been obtained by showing the
influence of the geometry on such problems. Particularly, the arising of
concentration phenomena in PDEs has motivated the
 development of
new methods in geometric analysis, see \cite{Au4}, \cite{DH} and \cite{H3} for
a complete survey. Our interest here is the study of the behavior of the
second Riemannian $L^p$-Sobolev best constant when the metric changes and some consequences such as existence and compactness results on extremal functions involving sets of Riemannian metrics.

Let $(M,g)$ be a smooth compact Riemannian manifold of dimension $n
 \geq 2$.
 For $1 \leq p < n$, we denote by $H^p_1(M)$ the standard first-order
 Sobolev space defined as the completion of $C^{\infty}(M)$ with
 respect to
 the norm

\[
||u||_{H^p_1(M)} = \left( \int_{M} |\nabla_g u|^p\; dv_g  + \int_{M}
 |u|^p\; dv_g \right)^{1/p},
\]

\n where $dv_g$ denotes the Riemannian volume element of $g$. The Sobolev embedding theorem ensures that the inclusion $H^p_1(M)
 \subset L^{p^*}(M)$ is continuous for $p^* = \frac{np}{n-p}$. Thus,
 there exist constants $A, B \in \R$ such that for any $ u \in
 H^p_1(M)$,

\begin{gather} \label{AB}
\left(\int_{M} |u|^{p^*}\; dv_g\right)^{p/p^*} \leq A \int_{M}
 |\nabla_g u|^p \; dv_g + B \int_{M} |u|^p\; dv_g\ . \tag{$I^p_g(A,B)$}
\end{gather}

The first $L^p$-Sobolev best constant associated to (\ref{AB}) is defined by

\[
A_0(p,g) = \inf \{ A \in \R:\; \mbox{ there exists} \hspace{0,18cm} B
 \in \R \hspace{0,18cm} \mbox{such that} \hspace{0,18cm} (I^p_g(A,B))
 \hspace{0,18cm} \mbox{is valid}\}
\]

\n and, by Aubin \cite{Au2}, its value is given by

\[
K(n,p)^p = \sup_{u \in {\cal D}^p_1(\R^n) \setminus \{0\}}
 \frac{\left(\int_{\R^n}
 |u|^{p^*}\; dx\right)^{p/p^*}}{\int_{\R^n} |\nabla u|^p \;
 dx}
\]

\n and ${\cal D}^p_1(\R^n)$ is the completion of $C^\infty_0(\R^n)$
under the norm

\[
||u||_{{\cal D}^p_1(\R^n)} = \left(\int_{\R^n} |\nabla u|^p\;
 dx\right)^{1/p}\ .
\]

\n In particular, the first best constant $A_0(p,g)$ does not depend on
 the geometry.

The first optimal Riemannain $L^p$-Sobolev inequality on $H^p_1(M)$ states that

\begin{gather} \label{A-opt}
\left(\int_{M} |u|^{p^*}\; dv_g\right)^{p/p^*} \leq K(n,p)^p \int_{M}
 |\nabla_g u|^p \; dv_g + B \int_{M} |u|^p\; dv_g
\tag{$I^p_{g,opt}$}
\end{gather}

\n for some constant $B \in \R$. The validity of (\ref{A-opt}) has been
 established by Hebey and Vaugon \cite{HV1} in the case $p=2$,
 independently, by Aubin and Li \cite{AuLi} and Druet \cite{D2} in the
 case $1 < p
 < 2$, and by Druet \cite{D4} in the case $p = 1$.

For $1 \leq p \leq 2$, define then the second $L^p$-Sobolev best
 constant by

\[
B_0(p,g) = \inf \{ B \in \R:\; (I^p_{g,opt}) \hspace{0,18cm} \mbox{is
 valid} \}\ .
\]

\n Clearly, for any $u \in H_1^p(M)$, one has

\begin{gather} \label{B-opt}
\left(\int_{M} |u|^{p^*}\; dv_g\right)^{p/p^*} \leq K(n,p)^p
\int_{M} |\nabla u|_g^p \; dv_g + B_0(p,g) \int_{M} |u|^p\; dv_g\ .
\tag{$J^p_{g,opt}$}
\end{gather}

\n Note that (\ref{B-opt}) is sharp with respect to both the first and
 second best constants in the sense that none of them can be lowered.
 The
 inequality (\ref{B-opt}) is called the second optimal Riemannain
 $L^p$-Sobolev inequality.

On the contrary of the first best constant, the second one depends strongly on the geometry. In fact, note that $B_0(p, \lambda g) = \lambda^{-1} B_0(p,g)$ for any constant $\lambda > 0$. An interesting remark is that the arguments used in the works \cite{AuLi}, \cite{DjDr}, \cite{D2}, \cite{D4} and \cite{HV1} rely only on the continuity of derivatives up to second order of the components of $g$. Thus, a natural question is to know if $B_0(p,g)$ depends continuously on the metric $g$ in the $C^2$-topology and if this topology is sharp.

Let $M$ be a smooth compact manifold of dimension $n \geq 2$. Denote by ${\cal M}_2$ the space of smooth Riemannian metrics on $M$ endowed with the $C^2$-topology and by ${\cal M}_\infty$ the space of smooth Riemannian metrics on $M$ endowed with the usual Fr\'{e}chet topology. We provide some answers for $1 \leq p \leq 2$ to the above question in the following theorems:

\begin{teor} \label{cont2}
Let $M$ be a smooth compact manifold of dimension $n$. If $n \geq 4$, then the map $g \in {\cal M}_2 \mapsto B_0(2,g)$ is continuous. Moreover, the $C^2$-topology is sharp.
\end{teor}

\begin{teor} \label{cont1}
Let $M$ be a smooth compact manifold of dimension $n$. If $n \geq 2$ and $1 \leq p < \min\{2,\sqrt{n}\}$, then the map $g \in {\cal M}_2 \mapsto B_0(p,g)$ is continuous.
\end{teor}

A direct consequence is:

\begin{cor} \label{frechet}
Let $M$ be a smooth compact manifold of dimension $n$. If either $1 \leq p \leq 2$ and $n \geq 4$ or $1 \leq p < \sqrt{n}$ and $n = 2, 3$, then the map $g \in {\cal M}_\infty \mapsto B_0(p,g)$ is continuous.
\end{cor}

The continuity question treated here is connected to the extremal functions
$C^0$-compactness and a uniformity problem as follows. An extremal function
 of (\ref{B-opt}) is a nonzero function $u_0 \in H^p_1(M)$ such that

\[
\left(\int_{M} |u_0|^{p^*}\; dv_g\right)^{p/p^*}= K(n,p)^p \int_{M}
 |\nabla_g u_0|^p \; dv_g + B_0(p,g) \int_{M} |u_0|^p\; dv_g\ .
\]

Let $G \subset {\cal M}_2$. Set $E_p(G) = \bigcup _{g\in G}E_p(g)$, where
 $E_p(g)$ denotes the set of
 all extremal functions of (\ref{B-opt}) with unit $L^{p^*}$-norm.

Consider a subset $G \subset {\cal M}_2$ such that

\[
B_0(2,g) > \frac{n-2}{4(n-1)}K(n,2)^2 \max_{M}Scal_{g}
\]

\n for all metric $g \in G$. By Theorem 1 of \cite{DjDr}, $E_2(g)$ is non-empty for
 all $g \in G$.

Theorem \ref{cont2} then implies the following compactness result:

\begin{cor} \label{compct2}
Let $n \geq 4$. If $G$ is compact in the $C^2$-topology, then $E_2(G)$ is compact in the $C^0$-topology.
\end{cor}

Let $G \subset {\cal M}_2$. If $n \geq 2$ and $1 \leq p <
 \min\{2,\sqrt{n}\}$, by Theorem 2 of \cite{DjDr}, $E_p(g)$ is non-empty and compact in the $C^0$-topology for all $g
 \in G$.

As a consequence of Theorem \ref{cont1}, we have:

\begin{cor} \label{compct1}
Let $n \geq 2$ and $1 \leq p < \min\{2,\sqrt{n}\}$. If $G$ is compact in the $C^2$-topology, then $E_p(G)$ is compact in the $C^0$-topology.
\end{cor}

Given a subset $G \subset {\cal M}_2$, the uniformity problem associated
 to (\ref{A-opt}) consists in knowing if there exists a constant $B>0$
 such that for any $u \in H^p_1(M)$ and any $g \in G$,

\begin{gather} \label{A-G-opt}
\left(\int_{M} |u|^{p^*}\; dv_g\right)^{p/p^*}\leq K(n,p)^p \int_{M}
 |\nabla_g u|^p \; dv_g + B \int_{M} |u|^p\; dv_g\ . \tag{$I^p_{g}(G)$}
\end{gather}

\n The existence of a such constant plays an important role in the
 study of Perelman's local non-collapsing properties along the Ricci
 flow.
 Recent advances in this direction have been obtained in \cite{Hsu},
 \cite{Ye} and \cite{Zhang}. In this context, $G$ represents the image
 of
 the flow in the space of metrics. The answer to this question clearly
 relies on properties of the set $G$. For example, as a consequence of
 Theorems \ref{cont2} and \ref{cont1}, if either $1 \leq p \leq 2$ and $n
 \geq
 4$ or $1 \leq p < \sqrt{n}$ and $n = 2, 3$, and $G$ is compact in the
 $C^2$-topology, then a such constant $B > 0$ exists. In this
 case, we define

\[
B_0(p,G) = \inf \{ B \in \R:\; (I^p_{g}(G)) \hspace{0,18cm} \mbox{is
 valid for all} \hspace{0,18cm} g \in G \}\ .
\]

\n Clearly,

\begin{gather} \label{B-G-opt}
\left(\int_{M} |u|^{p^*}\; dv_g\right)^{p/p^*}\leq K(n,p)^p \int_{M}
 |\nabla_g u|^p \; dv_g + B_0(p,G) \int_{M} |u|^p\; dv_g
 \tag{$I^p_{g,opt}(G)$}
\end{gather}

\n and

\[
B_0(p,G) = \sup_{g \in G} B_0(p,g)\ .
\]

\n Note that if (\ref{B-G-opt}) admits an extremal function for some metric $g
 \in G$, then $B_0(p,G) = B_0(p,g)$.

Existence results on extremal functions of (\ref{B-G-opt}) follow from results of \cite{DjDr} and from Theorems \ref{cont2} and \ref{cont1}.

Precisely, we have:

\begin{cor} \label{uniform2}
Let $n \geq 4$ and $G \subset {\cal M}_2$ be such that

\[
B_0(2,g) > \frac{n-2}{4(n-1)}K^2(n,2) \max_{M}Scal_{g}
\]

\n for all metric $g \in G$. If $G$ is compact in the $C^2$-topology, then ($I^2_{g,opt}(G)$) admits at least an extremal function.
\end{cor}

\begin{cor} \label{uniform1}
Let $n \geq 2$ and $1 \leq p < \min\{2,\sqrt{n}\}$. If $G$ is compact in
 the $C^2$-topology, then (\ref{B-G-opt}) admits at least an
 extremal
 function.
\end{cor}

The proofs of Theorems \ref{cont2} and \ref{cont1} are made by
 contradiction. If the conclusions fail, we naturally are
 led to two
 possible alternatives. One of them is directly eliminated according to
 the definition of second best constant. The other alternative implies
 the existence of minimizers, concentrating in a point, of functionals associated to a family of metrics. The idea then consists in
 performing a concentration refined study on these minimizers in order
 to
 obtain the second contradiction. The proofs are inspired in the works
 \cite{DjDr}, \cite{D2}, \cite{D4} and \cite{HV1}. New technical difficulties
 however
 arise when $g$ changes in ${\cal M}_2$. In all the study of
 concentration, we assume only $C^0$-convergence of metrics. The
 $C^2$-convergence is necessary only in the last
 step of the proofs. For $p = 2$, we construct a counter-example by showing that the $C^2$-topology is sharp
for the geometric continuity.\\

\section{Proof of Theorem \ref{cont2}}

\n  Consider initially a sequence $(g_\alpha)_\alpha \subset {\cal M}_2$
 converging to $g \in {\cal M}_2$ in the $C^0$-topology. The $C^2$-convergence will be used later
 in the last step of this proof. Suppose, by contradiction, that there
 exists $\varepsilon_0 >0$ such that $|B_0(2,g_{\alpha}) - B_0(2,g)| >
 \varepsilon _0$ for infinitely many $\alpha$. Then, at least, one of
 the situations holds:

\[
B_0(2,g) - B_0(2,g_{\alpha}) > \varepsilon _0 \ \ {\rm or}\ \
B_0(2,g_{\alpha}) - B_0(2,g) > \varepsilon _0
\]

\n for infinitely many $\alpha$. If the first alternative holds, then
 for any
 $u\in H^2_1(M)$,

\[
\left( \int_{M} |u|^{2^*}\; dv_{g_{\alpha}} \right)^{2/2^*}
 \leq
K(n,2)^2 \int_{M} |\nabla_{g_{\alpha}} u|^2 \; dv_{g_{\alpha}} +
 (B_0(2,g)-\epsilon _0) \int_{M} u^2\;
dv_{g_{\alpha}}\ .
\]

\n Letting $\alpha \rightarrow +\infty$ in this inequality, one finds

\[
\left( \int_{M} |u|^{2^*}\; dv_{g} \right)^{2/2^*} \leq K(n,2)^2
 \int_{M} |\nabla_{g} u|^2 \; dv_{g} + (B_0(2,g)-\varepsilon_0)
 \int_{M} u^2\;
 dv_{g},
\]

\n and this clearly contradicts the definition of $B_0(2,g)$. Suppose
 then that
 the second situation holds, i.e. $B_0(2,g) +
\varepsilon_0 < B_0(2,g_{\alpha})$ for infinitely many $\alpha$. For
 each $\alpha > 0$, consider the functional

\[
J_{\alpha }(u)= \int_{M} |\nabla_{g_{\alpha }} u|^2 \; dv_{g_{\alpha }}
 +(B_0(2,g) +\varepsilon _0) K(n,2)^{-2} \int_{M}  u^2 \;
 dv_{g_{\alpha}}
\]

\n defined on the set $\Lambda_{\alpha} =\{ u\in H^2_1(M):\;  \int_{M}
 |u|^{2^*} \; dv_{g_{\alpha }} = 1\}$. From the definition of
 $B_0(2,g_{\alpha})$, one has

\[
\lambda_{\alpha} := \inf _{u \in \Lambda_{\alpha}} J_{\alpha}(u) < K(n,2)^{-2}\
 .
\]

\n From this inequality, we find a nonnegative minimizer $u_{\alpha}
 \in
 \Lambda_{\alpha}$ for $\lambda_{\alpha}$. The Euler-Lagrange
 equation for $u_{\alpha}$ is

\begin{gather} \label{EL}
-\Delta _{g_{\alpha}}u_{\alpha} + (B_0(2,g) + \varepsilon_0)
 K(n,2)^{-2}u_{\alpha} = \lambda_{\alpha} u_{\alpha }^{2^*-1},
 \tag{$E_{\alpha}$}
\end{gather}

\n where $\Delta_{g_{\alpha}}u = {\rm
 div}_{g_{\alpha}}(\nabla_{g_{\alpha}} u)$ denotes the Laplace-Beltrami
 operator associated to the metric
 $g_{\alpha}$. From the classical elliptic theory, it follows then that
 $u_{\alpha} \in C^{\infty}(M)$ and $u_{\alpha }>0$ on $M$. Our goal
 now is
 to study the sequence $(u_{\alpha})_{\alpha}$ as
 $\alpha \rightarrow +\infty$. Note first that

\[
\int_{M} |\nabla_{g_{\alpha}} u_{\alpha}|^2 \; dv_{g_{\alpha}} +
 (B_0(2,g) + \varepsilon_0) K(n,2)^{-2} \int_{M} u_{\alpha}^2\;
 dv_{g_{\alpha}} = \lambda _{\alpha} < K(n,2)^{-2}
\]

\n and there exists a constant $c > 0$, independent of $\alpha$, such
 that

\[
\int_{M} u_{\alpha}^2\; dv_{g} \leq c \int_{M} u_{\alpha}^2\;
dv_{g_{\alpha}}
\]

\n and

\[
\int_{M} |\nabla_{g} u_{\alpha}|^2 \; dv_{g}  \leq c \int_{M}
 |\nabla_{g_{\alpha}} u_{\alpha}|^2 \; dv_{g_{\alpha}}
\]

\n for $\alpha > 0$ large. Clearly, these inequalities imply that
 $(u_{\alpha})_{\alpha}$ is bounded in $H^2_1(M)$. So, there exists $u
 \in H^2_1(M)$, $u\geq 0$, such that
 $u_{\alpha} \rightharpoonup u$ weakly in $H^2_1(M)$ and also
 $\lambda_{\alpha}
 \rightarrow \lambda$ as $\alpha \rightarrow + \infty $, with $0 \leq
 \lambda \leq K(n,2)^{-2}$, up to a subsequence. By the Sobolev embedding
 compactness theorem, we also have

\[
\int_{M} u_{\alpha}^q \; dv_{g_{\alpha}} \rightarrow  \int_{M} u^q \;
 dv_{g}
\]

\n for any $1 \leq q < 2^*$. Letting then $\alpha \rightarrow + \infty$
 in the equation (\ref{EL}) and using that $g_\alpha$ converges to $g$
 in $C^0$, we find that $u$ satisfies

\begin{gather} \label{ELLIM}
-\Delta_g u + (B_0(2,g) + \varepsilon_0) K(n,2)^{-2} u = \lambda
 u^{2^*-1}\ . \tag{$E$}
\end{gather}

\n If $u \not \equiv 0$, then ($J^2_{g,opt}$)
 and
 (\ref{ELLIM}) imply

\[
\left(\int_{M} u^{2^*}\; dv_g\right)^{2/2^*} < K(n,2)^2 \int_{M}
 |\nabla_g u|^2 \; dv_g + (B_0(2,g) + \varepsilon_0) \int_{M} u^2\;
 dv_g
\]

\[
= K(n,2)^2\lambda  \int_{M} u^{2^*}\; dv_g \leq \int_{M} u^{2^*}\;
 dv_g,
\]

\n so that $||u||_{2^*} > 1$. But, this is a contradiction, since

\[
\int _M u^{2^*}dv_{g} \leq \liminf_{\alpha \rightarrow + \infty} \int_M
 u_{\alpha}^{2^*}dv_{g_{\alpha}} = 1\ .
\]

\n We then assume that $u \equiv 0$ on $M$ and will obtain another
 contradiction. In this case, we claim that $\lambda_{\alpha}
 \rightarrow
 K(n,2)^{-2}$ as $\alpha \rightarrow + \infty$. In fact, using that
 $u_\alpha \in \Lambda_\alpha$ and $g_\alpha \rightarrow g$ in $C^0$, one
 gets

\[
\lim_{\alpha \rightarrow + \infty} \int_{M} u_\alpha^{2^*}\; dv_g = 1
\]

\n and

\[
\lim_{\alpha \rightarrow + \infty} \int_{M} u_{\alpha}^2\;
 dv_{g} = 0\ .
\]

\n Letting $\alpha \rightarrow +\infty$ in

\[
\left( \int_{M} u_\alpha^{2^*}\; dv_g \right)^{2/2^*} \leq K(n,2)^2
 \int_{M} |\nabla_g u_\alpha |^2 \; dv_g + B_0(2,g) \int_{M} u_\alpha^2\;
 dv_g
\]

\n and using the limits above, one finds

\[
\liminf_{\alpha \rightarrow +\infty} \int_{M} |\nabla_g u_\alpha |^2 \;
 dv_g \geq K(n,2)^{-2}\ .
\]

\n Clearly, the $C^0$-convergence of $g_\alpha$ then implies

\[
\liminf_{\alpha \rightarrow +\infty} \int_{M} |\nabla_{g_{\alpha}}
 u_\alpha |^2 \; dv_{g_{\alpha}} \geq K(n,2)^{-2}\ .
\]

\n The claim follows then from

\[
\limsup_{\alpha \rightarrow +\infty} \int_{M} |\nabla_{g_{\alpha}}
 u_\alpha |^2 \; dv_{g_{\alpha}} \leq
 \limsup_{\alpha \rightarrow +\infty} \lambda_{\alpha} \leq
 K(n,2)^{-2}\ .
\]

In the sequel, we divide the proof into six steps. Several possibly
 different positive constants, independent of $\alpha$, will be
 denoted by $c$.

We say that $x \in M$ is a point of concentration of $(u_{\alpha})_{\alpha}$ if,
 for any $\delta >0$,

\[
\limsup \limits_{\alpha \rightarrow +\infty} \int_{B_g (x, \delta)}
 u_{\alpha }^{2^*} \; dv_{g_{\alpha}} > 0 \ .
\]

\n {\bf Step 1:} The sequence $(u_{\alpha})_{\alpha}$ possesses exactly one
 point of concentration $x_0$, up to a subsequence.\\

\n {\bf Proof:} The existence of, at least, one point of concentration
 follows directly from the compactness of $M$, since $u_{\alpha} \in
 \Lambda_\alpha$. Conversely, let
 $x_0$ be a point of concentration of $(u_{\alpha})_{\alpha}$. Let $\delta > 0$
 small
 and consider a smooth function $\eta \in
 C^{\infty}_0(B_g(x_0,\delta))$
 such that $0 \leq \eta \leq 1$ and $\eta =1$ in $B_g(x_0,\delta /2)$.
 Multiplying (\ref{EL}) by $\eta^{2} u_{\alpha}^k$, $k > 1$, and
 integrating over $M$, one has

\begin{equation}\label{e1}
- \int_M \eta ^{2} u_{\alpha }^k \Delta_{g_\alpha} u_{\alpha} \;
 dv_{g_\alpha} + (B_0(2,g) + \varepsilon_0) K(n,2)^{-2} \int_M \eta^{2}
 u_{\alpha}^{k+1}\; dv_{g_\alpha} = \lambda_\alpha \int_M \eta^{2}
 u_{\alpha}^{k + 2^* - 1}\; dv_{g_\alpha}\ .
\end{equation}

\n For each $\varepsilon >0$, there exists a constant $c_{\varepsilon}
 > 0$, independent of $\alpha$, since $g_\alpha \rightarrow g$ in
 $C^0$, such that

\[
\int_M |\nabla_{g_\alpha} (\eta u_{\alpha}^{\frac{k+1}{2}})|^{2} \;
 dv_{g_\alpha} \leq \frac{(k+1)^2}{4} (1 + \varepsilon) \int_M \eta^{2}
 u_{\alpha}^{k-1} |\nabla_{g_\alpha} u_{\alpha}|^{2} \; dv_{g_\alpha}
\]

\[
+ c_{\varepsilon} ||\nabla \eta||_{\infty}^{2} \int_M u_{\alpha
 }^{k+1}\; dv_{g_\alpha}
\]

\n for $\alpha > 0$ large. By direct integration, we have

\[
- \int_M \eta^{2} u_{\alpha}^k \Delta_{g_\alpha} u_{\alpha}\;
 dv_{g_\alpha} \geq k \int_M \eta^{2} u_{\alpha }^{k-1}
 |\nabla_{g_\alpha}
 u_{\alpha}|^{2} \; dv_{g_\alpha} - \int_M u_{\alpha}^k
 |\nabla_{g_\alpha}
 u_{\alpha}| |\nabla_{g_\alpha} (\eta ^{2})| \; dv_{g_\alpha},
\]

\n so that, together with (\ref{e1}),

\begin{equation} \label{e2}
\int_M |\nabla_{g_\alpha} ( \eta u_{\alpha}^{\frac{k+1}{2}} )|^{2} \;
 dv_{g_\alpha} \leq  \frac{(k+1)^2}{4k} (1 + \varepsilon)
  \lambda_{\alpha} \int_M \eta^{2} u_{\alpha}^{k + 2^* -1} \;
 dv_{g_\alpha}
\end{equation}

\[
+ \frac{(k+1)^2}{4k} (1 + \varepsilon) \int_M u_{\alpha}^k
 |\nabla_{g_\alpha} u_{\alpha }| |\nabla_{g_\alpha} (\eta ^{2})|\;
 dv_{g_\alpha} +
 c_{\varepsilon} ||\nabla \eta||_{\infty}^{2} \int_M u_{\alpha}^{k+1}\;
 dv_{g_\alpha}\ .
\]

\n From the H\"{o}lder inequality, one has

\begin{equation} \label{e3}
\int_M \eta^{2} u_{\alpha}^{k+2^*-1} \; dv_{g_\alpha} \leq \left (
 \int_M ( \eta u_{\alpha}^{\frac{k+1}{2}} )^{2^*} \; dv_{g_\alpha}
 \right)^{2/2^*} \left( \int_{B_g(x_0,\delta)} u_{\alpha}^{2^*} \;
 dv_{g_\alpha}
 \right)^{1-2/2^*}
\end{equation}

\n and

\begin{equation} \label{e4}
\int_M u_{\alpha}^k |\nabla_{g_\alpha} u_{\alpha}| |\nabla_{g_\alpha}
 (\eta ^{2})| \; dv_{g_\alpha} \leq 2 ||\nabla \eta||_{\infty} \left(
 \int_M |\nabla_{g_\alpha} u_{\alpha}|^{2} \; dv_{g_\alpha}
 \right)^{1/2}
 \left( \int_M u_{\alpha}^{2k} \; dv_{g_\alpha} \right)^{1/2} \ .
\end{equation}

\n For each $\varepsilon >0$, there exists a constant $d_{\varepsilon}
 > 0$, independent of $\alpha$, such that

\begin{equation} \label{e5}
\left( \int_M ( \eta u_{\alpha}^{\frac{k+1}{2}} )^{2^*} \;
 dv_{g_\alpha} \right)^{2/2^*} \leq (K(n,2)^{2} + \varepsilon) \int_M
 |\nabla_{g_\alpha} ( \eta u_{\alpha}^{\frac{k+1}{2}} )|^{2} \;
 dv_{g_\alpha} +
 d_{\varepsilon} \int_M u_{\alpha}^{k+1} \; dv_{g_\alpha}
\end{equation}

\n for $\alpha > 0$ large. Here is used that $(1 - \varepsilon) g \leq
 g_\alpha \leq (1 + \varepsilon) g$ in the bilinear forms sense. From $J_\alpha(u_\alpha) < K(n,2)^{-2}$, one has

\begin{equation} \label{e6}
\left( \int_M |\nabla_{g_\alpha} u_{\alpha }|^{2} \; dv_{g_\alpha}
  \right)^{1/2} \leq \left( \lambda_{\alpha} \int_M u_{\alpha}^{2^*} \;
 dv_{g_\alpha} \right)^{1/2} \leq K(n,2)^{-1}\ .
\end{equation}

\n So, putting together (\ref{e2}), (\ref{e3}), (\ref{e4}), (\ref{e5})
 and (\ref{e6}), one finds

\begin{equation}\label{e7}
A_{\alpha} \left( \int_M ( \eta u_{\alpha}^{\frac{k+1}{2}} )^{2^*} \;
 dv_{g_\alpha} \right)^{2/2^*} \leq B \int_M u_{\alpha}^{k+1} \;
 dv_{g_\alpha} + C \left( \int_M u_{\alpha}^{2 k} \; dv_{g_\alpha}
 \right)^{1/2},
\end{equation}

\n where

\[
A_{\alpha} = 1 - \frac{(k+1)^2}{4k} (1 + \varepsilon)^2
 \lambda_{\alpha} K(n,2)^{2} \left( \int_{B_g(x_0,\delta)} u_{\alpha
 }^{2^*} \;
 dv_{g_\alpha} \right)^{1-2/2^*},
\]

\[
B = K(n,2)^{2} (1 + \varepsilon) c_{\varepsilon} ||\nabla
 \eta||_{\infty}^{2} + d_{\varepsilon}
\]

\n and

\[
C = 2 \frac{(k+1)^2}{4k} (1 + \varepsilon)^2 ||\nabla \eta||_{\infty}
 K(n,2)\ .
\]

\n Since $x_0$ is a point of concentration of $(u_{\alpha})_{\alpha}$, we have

\[
\limsup \limits_{\alpha \rightarrow +\infty} \left(
 \int_{B_g(x_0,\delta)} u_{\alpha}^{2^*}\; dv_{g_\alpha}
 \right)^{1-2/2^*} = a > 0,
\]

\n with $a \leq 1$, since $u_{\alpha} \in \Lambda_\alpha$. We
 claim that $a=1$ for all $\delta > 0$. In fact, if $a<1$ for some
 $\delta>0$, taking $\varepsilon > 0$ small enough and $k > 1$ close to
 $1$
 such that $A_{\alpha} > A$, where $A$ is a positive constant and
 independent
 of $\alpha$. Since the right-hand side of (\ref{e7}) is bounded for
 $k$
 close to $1$, we find a constant $c > 0$, independent of $\alpha$,
 such that

\[
\left( \int_M ( \eta u_{\alpha}^{\frac{k+1}{2}} )^{2^*} \;
 dv_{g_\alpha} \right)^{2/2^*} \leq c
\]

\n for $\alpha > 0$ large. From the H\"{o}lder inequality, one has

\[
\int_{B_g(x_0,\frac{\delta}{2})} u_{\alpha }^{2^*} \; dv_{g_\alpha} =
 \int_{B_g(x_0,\frac{\delta}{2})} u_{\alpha}^{k+1}\; u_{\alpha}^{2^* -
 1
 - k} \; dv_{g_\alpha}
\]

\[
\leq \left( \int_M ( \eta u_{\alpha}^{\frac{k+1}{2}} )^{2^*} \;
 dv_{g_\alpha}
\right)^{2/2^*} \left( \int_M u_{\alpha}^{2^* - \frac{2^*(k-1)}{2^*-2}}
 \; dv_{g_\alpha}\right )^{1-2/2^*} \leq c \left( \int_M
 u_{\alpha}^{2^*-\frac{2^*(k-1)}{2^*-2}}\; dv_{g_\alpha}\right
 )^{1-2/2^*}\ .
\]

\n Choose $k$ close to $1$ such that $2 < 2^* -\frac{2^*(k-1)}{2^*-2}
 < 2^*$. Since $||u_{\alpha}||_{2} \rightarrow 0$, it follows then
 from an interpolation argument that

\[
\limsup \limits_{\alpha \rightarrow +\infty}
 \int_{B_g(x_0,\frac{\delta}{2})} u_{\alpha }^{2^*}\; dv_{g_\alpha} =
 0\ .
\]

\n But this clearly contradicts the fact that $x_0$ is a point of
 concentration. Therefore, $a=1$ and

\[
\limsup \limits_{\alpha \rightarrow +\infty} \int_{B_g(x_0,\delta)}
 u_{\alpha }^{2^*} \; dv_{g_\alpha} = 1
\]

\n for all $\delta > 0$. Since $u_{\alpha} \in \Lambda_\alpha$, it
 follows then that
$(u_{\alpha })_{\alpha}$ has exactly one point of
 concentration, up to a subsequence.\bl\\

\n {\bf Step 2:} Let $x_0 \in M$ be the unique point of concentration
 of
 $(u_{\alpha})_{\alpha}$. Then,

\begin{equation}\label{e9}
\lim \limits _{\alpha \rightarrow +\infty} u_{\alpha} = 0\ \ {\rm in}\
 \ C_{loc}^0(M \setminus \{x_0\})\ .
\end{equation}

\n {\bf Proof:} From (\ref{e7}), given $\overline{\Omega} \subset M
 \setminus \{x_0\}$, there exist constants $\varepsilon, c_1 > 0$,
 independent of $\alpha$, such that

\[
\int_{\Omega} u_{\alpha }^{2^*(1 + \varepsilon)} \; dv_{g_\alpha} \leq
 c_1
\]

\n for $\alpha > 0$ large. On the other hand, from the
 $C^0$-convergence of $g_\alpha$, we find constants
 $\gamma$ and $c_0$ such that $g_{\alpha} \geq \gamma \xi$, in the
 bilinear
 forms sense, and $||(g_{\alpha})_{ij}||_{C^0} \leq c_0$ for
 $\alpha > 0$ large, where $\xi $ stands for the Euclidean metric on
 $\R^n$.
 Finally, the conclusion (\ref{e9}) follows from a De Giorgi-Nash-Moser
 iterative scheme applied to (\ref{EL}). Here, it is important to note
 that the involved constants in this scheme
 depend only on $\gamma$, $c_0$ and $c_1$. We refer for instance to
 Serrin \cite{Se} for more details.\bl\\

Let $x_\alpha \in M$ be a maximum point of $u_\alpha$, i.e.
 $u_\alpha(x_\alpha) = ||u_\alpha||_{\infty}$. By the steps $1$ and
 $2$, one has $x_{\alpha} \rightarrow x_0$ as $\alpha \rightarrow
+\infty$.\\

\n {\bf Step 3:} For each $R > 0$, we have

\begin{equation} \label{first}
\lim _{\alpha \rightarrow  + \infty} \int _{B_{g_{\alpha}}(x_{\alpha},
 R \mu_{\alpha})} u_{\alpha }^{2^*}\; dv_{g_{\alpha }} = 1 -
 \varepsilon_{R}
\end{equation}

\n where $\mu_\alpha = ||u_\alpha||_{\infty}^{- 2^*/n}$ and
 $\varepsilon = \varepsilon_R \rightarrow 0$ as $R \rightarrow +
 \infty$.\\

\n {\bf Proof:} From

\[
1 = \int_{M} u_{\alpha}^{2^*}\; dv_{g_{\alpha }}\leq ||u_{\alpha
 }||_{\infty}^{2^* - 2} \int_{M} u_{\alpha}^2\; dv_{g_{\alpha}},
\]

\n we find $||u_{\alpha }||_{\infty} \rightarrow +\infty$ as $\alpha
 \rightarrow  +\infty $, since $ \int_{M} u_{\alpha }^2\; dv_{g_\alpha}
 \rightarrow
 0$. So, $\mu_{\alpha} \rightarrow 0$ as $\alpha \rightarrow
  + \infty$. Let $\exp_{x_\alpha}$ be the exponential map at
 $x_{\alpha}$ with respect to the metric $g$. Since $x_\alpha
 \rightarrow
 x_0$, there exists $\delta > 0$, independent of $\alpha$, such that
 $\exp_{x_{\alpha}}$ map $B(0,\delta) \subset \R^n$ onto
 $B_{g}(x_{\alpha},
 \delta)$ for $\alpha > 0$ large. For each $x\in B(0, \delta
 \mu_{\alpha}^{-1})$, set

\[
\tilde{g}_{\alpha}(x) = (\exp_{x_{\alpha}}^*g_{\alpha})(\mu_{\alpha}x)
\]

\n and

\[
\varphi_{\alpha}(x) = \mu_{\alpha}^{n/2^*} u_{\alpha}
 (\exp_{x_{\alpha}}(\mu_{\alpha}x))\ .
\]

\n As one easily checks,

\begin{gather} \label{ELT}
- \Delta_{\tilde{g}_{\alpha}} \varphi_{\alpha} + (B_0(2,g) +
 \varepsilon_0) K(n,2)^{-2} \mu_{\alpha}^2 \varphi_{\alpha} =
 \lambda_{\alpha}
 \varphi_{\alpha}^{2^* - 1}\ .  \tag{$\tilde{E}_{\alpha}$}
\end{gather}

\n Clearly,

\begin{equation} \label{9}
\tilde{g}_{\alpha} \rightarrow \xi\ \ {\rm in}\ \ C^0_{loc}(\R^n)\ .
\end{equation}

\n In particular, for each bounded open $\Omega \subset \R^n$, there
 exist constants
 $\gamma, c_0 > 0$ such that

\begin{equation} \label{10}
\tilde{g}_{\alpha} \geq \gamma \xi\ \ {\rm in}\ \ \Omega,
\end{equation}

\n in the bilinear forms sense, and

\begin{equation} \label{11}
||(\tilde{g}_{\alpha})_{ij}||_{C^0(\Omega)} \leq c_0
\end{equation}

\n for $\alpha > 0$ large. So, from (\ref{10}), there exists a constant
 $c > 0$ such that

\[
\int_{\Omega} |\nabla \varphi_{\alpha}|^2\; dv_{\xi} \leq c
 \int_{B(0,\delta \mu_{\alpha}^{-1})} |\nabla_{\tilde{g}_{\alpha}}
 \varphi_{\alpha}|^2\; dv_{\tilde{g}_{\alpha}} = c \int_{B(x_\alpha
 ,\delta)}
 |\nabla_{g_\alpha} u_{\alpha}|^2 \; dv_{g_\alpha} \leq  c K(n,2)^{-2}
\]

\n and

\[
\int_{\Omega} \varphi_{\alpha}^{2^*}\; dv_{\xi} \leq c \int_{B(0,\delta
 \mu_{\alpha}^{-1})} \varphi_{\alpha}^{2^*}\; dv_{\tilde{g}_{\alpha}} =
 c \int_{B(x_\alpha ,\delta)} u_{\alpha}^{2^*} \; dv_{g_\alpha} \leq
  c\ .
\]

\n Therefore, the sequence $(\varphi_{\alpha})_{\alpha}$, with
 $\alpha > 0$ large, is bounded in $H^2_1(\Omega)$ for any bounded open
 $\Omega
 \subset \R^n$, so that $\varphi_{\alpha} \rightharpoonup \varphi$
 weakly in
 $H^2_1(\Omega)$, $\varphi \geq
 0$, and $\int_{\Omega} \varphi_{\alpha}^q \; dv_{\xi}
 \rightarrow
  \int_{\Omega} \varphi^q \; dv_{\xi}$ for any $1 \leq q < 2^*$, up to
 a subsequence. Then, letting $\alpha \rightarrow + \infty$ in
 (\ref{ELT}),
 using (\ref{9}), $\lambda_\alpha \rightarrow K(n,2)^{-2}$ and
 $\mu_{\alpha} \rightarrow 0$, we conclude that $\varphi$ satisfies in
 the
 weak sense,

\begin{equation} \label{12}
- \Delta \varphi = K(n,2)^{-2} \varphi^{2^*-1}\ \ {\rm in}\ \ \R^n\ .
\end{equation}

\n Note also that $\varphi \in {\cal D}^2_1(\R^n)$. This last fact
 follows directly from

\[
\int_{\Omega} |\nabla \varphi_{\alpha}|^2\; dv_{\xi} \leq c
 K(n,2)^{-2}
\]

\n and $\varphi_{\alpha} \rightharpoonup \varphi$ in $H^2_1(\Omega)$.
 Thanks to (\ref{10}), (\ref{11}) and the bound
 of $(\varphi_{\alpha})_{\alpha}$ and $(\mu_{\alpha})_{\alpha}$, classical H\"{o}lder
 estimates on
 elliptic PDEs weak solutions (see \cite{LU}) can be applied to
 (\ref{ELT}), so that $(\varphi_{\alpha})_{\alpha}$ is uniformly bounded in
 $C^\beta(\overline{\Omega})$ for any
 bounded open $\Omega \subset \R^n$ and $\alpha > 0$ large. Therefore,
 $\varphi_{\alpha} \rightarrow \varphi$ in $C^0_{loc}(\R^n)$, up to a
 subsequence, so that $\varphi \not \equiv 0$, since $\varphi_{\alpha}
 (0) = 1$
 for all $\alpha$. From the equation (\ref{12}), one has

\[
\int_{\R^n} |\nabla \varphi|^2\; dv_{\xi} = K(n,2)^{-2} \int_{\R^n}
 \varphi^{2^*}\; dv_{\xi},
\]

\n since $\varphi \in {\cal D}^2_1(\R^n)$. So,

\[
K(n,2)^{-2} \left( \int_{\R^n} \varphi^{2^*}\; dv_{\xi} \right)^{2/2^*}
 \leq \int_{\R^n} |\nabla \varphi|^2\; dv_{\xi} = K(n,2)^{-2}
 \int_{\R^n} \varphi^{2^*} \; dv_{\xi},
\]

\n so that

\[
\int_{\R^n} \varphi^{2^*}\; dv_{\xi} \geq 1\ .
\]

\n On the other hand, since

\[
\int_{\Omega} \varphi_{\alpha}^{2^*}\; dv_{\tilde{g}_{\alpha}} \leq
 \int_{B(0, \delta \mu_{\alpha}^{-1})} \varphi_{\alpha}^{2^*}\;
 dv_{\tilde{g}_{\alpha}}  = \int_{B_{g_{\alpha}}(x_{\alpha}, \delta)}
 u_{\alpha}^{2^*} \; dv_{g_{\alpha}} \leq 1,
\]

\n we find $\int_{\R^n} \varphi^{2^*} \; dv_{\xi} = 1$, so that the
 conclusion of this step follows from the convergence

\[
\int_{B_{g_{\alpha}}(x_{\alpha }, R\mu_{\alpha})} u_{\alpha}^{2^*}\;
 dv_{g_{\alpha}} = \int_{B(0, R)} \varphi_{\alpha}^{2^*}\;
 dv_{\tilde{g}_{\alpha}} \rightarrow \int_{B(0, R)} \varphi^{2^*} \;
 dv_{\xi}\ .\ \ \ \bl
\]

\n {\bf Step 4:} There exists a constant $c > 0$, independent of
 $\alpha$, such that

\[
d_{g}(x,x_{\alpha})^{n/2^*} u_{\alpha}(x) \leq c
\]

\n for all $x \in M$ and $\alpha$ large, where $d_{g}$ stands
 for the distance with respect to the metric $g$.\\

\n {\bf Proof:} Set $\omega_{\alpha}(x) = d_{g}(x,x_{\alpha})^{n/2^*}
 u_{\alpha}(x)$ for $x \in M$ and suppose, by contradiction, that the
 conclusion of this step fails. In this case,

\[
\lim_{\alpha \rightarrow +\infty}||\omega_{\alpha }|| _{\infty} =
 +\infty,
\]

\n up to a subsequence. We next will derive a contradiction. Let
 $y_{\alpha} \in M$ be a maximum point of $\omega_{\alpha}$. Note that
 $u_{\alpha}(y_{\alpha}) \rightarrow +\infty$ and

\begin{equation} \label{1}
\lim_{\alpha \rightarrow + \infty}
 \frac{d_{g}(y_{\alpha},x_{\alpha})}{\mu_{\alpha}} = +\infty,
\end{equation}

\n since

\[
\frac{d_{g}(y_{\alpha},x_{\alpha})}{\mu _{\alpha}} =
 \frac{\omega_{\alpha }(y_{\alpha})^{2^*/n}}{\mu_{\alpha}
u_{\alpha}(y_{\alpha})^{2^*/n}} \geq
 \omega_{\alpha}(y_{\alpha})^{2^*/n} \ .
\]

\n Let $\exp_{x_\alpha}$ be the exponential map at
 $x_{\alpha}$ with respect to the metric $g$. For $x \in B(0,2)$, we define

\[
\hat{g}_{\alpha}(x) =
 (\exp_{y_{\alpha}}^{*}g_{\alpha})(u_{\alpha}(y_{\alpha})^{-2^*/n}x)\ .
\]

\n and

\[
v_{\alpha}(x) = u_{\alpha}(y_{\alpha})^{-1}
 u_{\alpha}(\exp_{y_{\alpha}}(u_{\alpha}(y_{\alpha})^{-2^*/n}x))
\]

\n We claim that the sequence $(v_{\alpha})_{\alpha}$ is uniformly bounded on
 $B(0,2)$ for $\alpha > 0$ large. In fact, there exists a constant $c >
 0$, independent of $\alpha$, such that for any $x \in B(0,2)$ and
 $\alpha > 0$ large, one has

\[
d_{g} (x_{\alpha}, \exp_{y_{\alpha}}(u_{\alpha}(y_{\alpha})^{-2^*/n}x))
 \geq d_{g}(x_{\alpha}, y_{\alpha}) - d_{g} (y_{\alpha},
 \exp_{y_{\alpha}}(u_{\alpha}(y_{\alpha})^{-2^*/n}x))
\]

\[
= d_{g}(x_{\alpha}, y_{\alpha}) - 2 u_{\alpha}(y_{\alpha})^{-2^*/n} =
 (1 - 2 \omega_\alpha(y_\alpha)^{-2^*/n} ) d_{g}(x_\alpha, y_\alpha)\ .
\]

\n Since $\omega_{\alpha}(y_{\alpha}) \rightarrow +\infty $ as $\alpha
 \rightarrow +\infty$, for $\alpha > 0$ large, one has

\begin{equation} \label{3}
d_{g} (x_{\alpha}, \exp_{y_{\alpha}}(u_{\alpha}(y_{\alpha})^{-2^*/n}x))
\geq \frac{1}{2}d_{g}(x_{\alpha}, y_{\alpha})\ .
\end{equation}

\n Hence,

\[
v_{\alpha}(x) = u_{\alpha}(y_{\alpha})^{-1}
 u_{\alpha}(\exp_{y_{\alpha}}(u_{\alpha}(y_{\alpha})^{-2^*/n}x))
\]

\[
\leq 2^{n/2^*} d_{g}(x_{\alpha},
 y_{\alpha})^{-n/2^*}u_{\alpha}(y_{\alpha})^{-1}
 \omega_{\alpha}(\exp_{y_{\alpha}}(u_{\alpha}(y_{\alpha})^{-2^*/n}x))
\]

\[
\leq 2^{n/2^*} d_{g}(x_{\alpha},
 y_{\alpha})^{-n/2^*}u_{\alpha}(y_{\alpha})^{-1}
 \omega_{\alpha}(y_{\alpha}) = 2^{n/2^*},
\]

\n so that

\begin{equation} \label{4}
||v_{\alpha}||_{L^{\infty}(B(0,2))} \leq 2^{n/2^*}\ .
\end{equation}

\n On the other hand, $v_{\alpha}$ satisfies

\[
-\Delta_{\hat{g}_{\alpha}} v_{\alpha} + B_{\alpha} v_{\alpha} =
 \lambda_{\alpha} v_{\alpha}^{2^*-1}\ \ {\rm in}\ \ B(0,2)
\]

\n for a certain constant $B_{\alpha} > 0$, so that

\begin{equation} \label{2}
-\Delta_{\hat{g}_{\alpha}} v_{\alpha} \leq \lambda_{\alpha}
 v_{\alpha}^{2^*-1}\ \ {\rm in}\ \ B(0,2)\ .
\end{equation}

\n Note also that

\[
\hat{g}_{\alpha} \rightarrow \xi\ \ {\rm in}\ \ C^0(\overline{B(0,2)}),
\]

\n so that there exist constants $\gamma, c_0 > 0$ such that

\begin{equation} \label{13}
\hat{g}_{\alpha} \geq \gamma \xi\ \ {\rm in}\ \ \Omega,
\end{equation}

\n in the bilinear forms sense, and

\begin{equation} \label{14}
||(\hat{g}_{\alpha})_{ij}||_{C^0(\overline{B(0,2)})} \leq c_0
\end{equation}

\n for $\alpha > 0$ large. Thanks to (\ref{4}), (\ref{13}) and (\ref{14}),
 the De Giorgi-Nash-Moser iterative scheme can be applied to (\ref{2}),
 so that

\[
v_{\alpha}(0) \leq \sup_{B(0,1)} v_{\alpha}(x) \leq  c \int_{B(0,2)}
 v_{\alpha}^{2^*} \; dv_{\hat{g}_{\alpha}} = c
 \int_{B_{g_{\alpha}}(y_{\alpha},2u_{\alpha}(y_{\alpha})^{-2^*/n})}
 u_{\alpha}^{2^*}\;
 dv_{g_{\alpha}}
\]

\n for some constant $c>0$ depending only on $\gamma$ and $c_0$. Since
 $v_{\alpha}(0) = 1$, the desired contradiction is then obtained by
 showing that the right-hand side integral converges to $0$ as $\alpha
 \rightarrow +\infty$. By the step $3$, it is sufficient then to show
 that

\[
B_{g_{\alpha}} (y_{\alpha}, 2u_{\alpha}(y_{\alpha})^{-2^*/n} ) \cap
 B_{g_{\alpha}} ( x_{\alpha}, R\mu _{\alpha} ) = \emptyset\ .
\]

\n Since $g_\alpha \rightarrow g$ in $C^0$ and $M$ is compact, there
 exists a constant $c > 0$, independent of $\alpha$, such that
 $d_{g_\alpha} \geq c d_{g}$ for $\alpha > 0$ large. Then, the assertion
 above follows
 directly from

\[
d_{g_\alpha} (x_\alpha, y_\alpha) u_{\alpha}(y_{\alpha})^{2^*/n} \geq
c d_{g} (x_\alpha, y_\alpha) u_{\alpha}(y_{\alpha})^{2^*/n}
\]

\[
= c w_\alpha (y_\alpha)^{2^*/n} \geq 2 + R u_{\alpha}(y_{\alpha})^{2^*/n}
\mu_\alpha = 2 + R u_{\alpha}(y_{\alpha})^{2^*/n}
 ||u_\alpha||_{\infty}^{- 2^*/n},
\]

\n which clearly holds for $\alpha > 0$ large, since $w_\alpha (y_\alpha)
 \rightarrow +\infty$.\bl\\

\n {\bf Step 5:} For each $\delta > 0$ small, one has

\begin{equation} \label{concent}
\lim_{\alpha \rightarrow +\infty} \frac{\int_{M \setminus
 B_{g}(x_0,\delta)} u_{\alpha}^{2}\; dv_{g}}{\int_{M} u_{\alpha}^2 \;
 dv_{g}} = 0\ .
\end{equation}

\n {\bf Proof:} First, by H\"{o}lder's inequality,

\[
\int_{M \setminus B_{g}(x_0, \delta)} u_{\alpha}^{2}\; dv_{g} \leq
 \sup_{M \setminus B_{g}(x_0, \delta)} u_{\alpha}\; \int_M u_{\alpha}
 \;
 dv_{g}  \leq v_{g}(M)^{1/2} \sup_{M \setminus B_{g}(x_0, \delta)}
 u_{\alpha}\; \left(  \int_M u_{\alpha}^2 \; dv_{g}\right)^{1/2}\ .
\]

\n By the step 2, the $C^0$-convergence of $g_\alpha$,
 (\ref{10}) and (\ref{11}), we can perform a De Giorgi-Nash-Moser
 iterative scheme in (\ref{EL}) and find a constant $c_1, c_2 > 0$,
 depending
 only on $\gamma$, $c_0$ and $\delta$, such that

\[
\sup_{M \setminus B_{g}(x_0, \delta)} u_{\alpha} \leq c_1 \int_M
 u_{\alpha} \; dv_{g_{\alpha}} \leq c_2 \int_M u_{\alpha} \; dv_{g}
\]

\n for $\alpha > 0$ large. From two inequalities above and (\ref{EL}), one
 finds

\begin{equation} \label{6}
\int_{M \setminus B_{g}(x_0, \delta)} u_{\alpha}^{2}\; dv_{g} \leq c
 \left(  \int_M u_{\alpha}^2 \; dv_{g}\right)^{1/2} \int_M
 u_{\alpha}^{2^*
 - 1} \; dv_{g}
\end{equation}

\n for $\alpha$ large. Now we analyze two situations. If $n = 4$, then

\[
\frac{ \int_{M \setminus B_{g}(x_0, \delta)} u_{\alpha}^{2} \;
 dv_{g}}{\int_{M} u_{\alpha}^{2} \; dv_{g}} \leq c ||u_{\alpha}||_{2}
 \rightarrow 0,
\]

\n since $2^* - 1 = 2$. Else, if $n > 4$, then $2^* - 1 > 2$. In this
 case, applying a H\"{o}lder type inequality and using that $u_{\alpha}
 \in \Lambda_\alpha$, one arrives at

\[
\frac{\int_{M \setminus B_{g}(x_0, \delta)} u_{\alpha}^{2} \;
 dv_{g}}{\int_{M} u_{\alpha}^{2} \; dv_{g}} \leq c
 ||u_{\alpha}||_{2}^{(n - 4)/2}
 \rightarrow 0\ .\ \ \ \bl
\]

\n {\bf Step 6:} Here is the final argument. Assume that
$g_\alpha$ converges to $g$ in the $C^2$-topology. Thus, we
have

\[
\liminf \limits _{\alpha \rightarrow +\infty}
{\rm inj}_{g_{\alpha}}(M)>0,
\]

\n where ${\rm inj}_{g_{\alpha}}(M)$ denotes the injectivity radius of
$(M, g_{\alpha})$. So, there exists $\delta > 0$ small enough, independent of $\alpha$, such that
$B_{g_{\alpha}}(x_{\alpha},\delta)$ is a geodesic ball for all
$\alpha > 0$ large. Moreover, if $\exp_{x_{\alpha}, g_\alpha}$ denote the
exponential map at $x_{\alpha}$ with respect to the metric
$g_{\alpha}$, then $\exp_{x_{\alpha}, g_\alpha} \circ \exp_{x_0, g}^{-1}$ converges to the identity map $id: B(0, \delta) \rightarrow \R^n$ in the $C^3$-topology. For each $x\in B(0,\delta)$, we set

\[
h_{\alpha}(x)= \exp_{x_{\alpha}, g_\alpha}^*g_{\alpha}(x)
\]

\n and

\[
v_{\alpha}(x)=u_{\alpha} ( \exp_{x_{\alpha}, g_\alpha}(x) )\ .
\]

\n Let $\eta \in C^{\infty }_0(B(0,\delta ))$ be such that $\eta =1$ on
$B(0,\frac{\delta}{2} )$ and $|\nabla \eta|\leq c\delta ^{-1}$. In the sequel, $c$ denotes a positive constant, independent of $\alpha$
and $\delta$. From the Euclidean $L^2$-Sobolev inequality, one has

\begin{equation}\label{new1}
\left (\int _{B(0,\delta )}\left ( \eta v_{\alpha}  \right
)^{2^*}\; dv_{\xi} \right )^{2/2^*}\leq K(n,2)^2\int
_{B(0,\delta )}|\nabla \left ( \eta v_{\alpha}  \right
)|^{2}\; dv_{\xi}\ .
\end{equation}

\n As easily one checks,

\[
\int_{B(0,\delta )}|\nabla ( \eta v_{\alpha} )|^{2}\; dv_{\xi}\leq \int _{B(0,\delta )}\eta^2 v_{\alpha} \Delta v_{\alpha}\; dv_{\xi} + c \delta^{-2}\int _{B(0,\delta) \setminus
B(0,\frac{\delta }{2})}v_{\alpha}^2\; dv_{\xi}\,.
\]

\n We also have

\[
-\Delta v_{\alpha} = -\Delta
_{h_{\alpha}}v_{\alpha}+(h_{\alpha}^{ij}-\delta _{ij})\partial
_{ij}v_{\alpha}-h_{\alpha}^{ij}\Gamma (h_{\alpha})^k_{ij}\partial
_{k}v_{\alpha}\,,
\]

\n where $\delta _{ij}$ is the Kronecker's symbol and $\Gamma
(h_{\alpha})^k_{ij}$ denotes the Christoffel's symbols of the
Levi-Civita's connection associated to the metric $h_{\alpha}$. This
gives

\[
\int_{B(0,\delta )}|\nabla  ( \eta v_{\alpha} )|^{2}\; dv_{\xi} \leq -\int _{B(0,\delta )}\eta^2 v_{\alpha}\Delta
_{h_{\alpha}}v_{\alpha}\; dv_{\xi}+ c \delta ^{-2} \int
_{B(0,\delta)\setminus B(0,\frac{\delta }{2})} v_{\alpha}^2\; dv_{\xi}
\]

\[
+\int _{B(0,\delta)}\eta ^2v_{\alpha}(h_{\alpha}^{ij}-\delta
_{ij})\partial _{ij}v_{\alpha}\; dv_{\xi}-\int _{B(0,\delta)}\eta
^2v_{\alpha}h_{\alpha}^{ij}\Gamma (h_{\alpha})^k_{ij}\partial
_{k}v_{\alpha}\; dv_{\xi},
\]

\n so that integrating by parts, using (\ref{EL}) and $\lambda _{\alpha}<K(n,2)^{-2}$, we get

\[
\int _{B(0,\delta )}|\nabla ( \eta v_{\alpha}
)|^{2}\; dv_{\xi} \leq K(n,2)^{-2}\int _{B(0,\delta)}\eta
^2v_{\alpha}^{2^*}\; dv_{\xi}-(B_0(2,g)+ \varepsilon _0)K(n,2)^{-2}\int
_{B(0,\delta)} (\eta v_{\alpha} )^2\; dv_{\xi}
\]

\[
+c\delta ^{-2}\int _{B(0,\delta)\setminus B(0,\frac{\delta
}{2})}v_{\alpha}^2\; dv_{\xi}- \int _{B(0,\delta)}\eta
^2(h_{\alpha}^{ij}-\delta _{ij})\partial _{i}v_{\alpha}\partial
_{j}v_{\alpha}\; dv_{\xi}
\]

\[
+\frac{1}{2}\int _{B(0,\delta)} (
\partial _k h_{\alpha}^{ij}\Gamma (h_{\alpha})^k_{ij} +\partial _{ij}h_{\alpha}^{ij} )
(\eta v_{\alpha})^2\; dv_{\xi}\ .
\]

\n From (\ref{new1}), one then obtains

\[
(B_0(2,g) + \varepsilon _0)\int _{B(0,\delta)} (\eta v_{\alpha}
 )^2\; dv_{\xi}\leq \int _{B(0,\delta)}\eta
^2v_{\alpha}^{2^*}\; dv_{\xi}-\left ( \int _{B(0,\delta)} (\eta
v_{\alpha} )^{2^*}\; dv_{\xi} \right )^{2/2^*}
\]

\[
+c\delta ^{-2}\int _{B(0,\delta)\setminus B(0,\frac{\delta
}{2})}v_{\alpha}^2\; dv_{\xi}+\frac{K(n,2)^2}{2}\int
_{B(0,\delta)} \partial _k (h_{\alpha}^{ij}\Gamma (h_{\alpha})^k_{ij} + \partial _{ij}h_{\alpha}^{ij} )
(\eta v_{\alpha})^2\; dv_{\xi}
\]

\[
-K(n,2)^2 \int _{B(0,\delta)}\eta ^2(h_{\alpha}^{ij}-\delta
_{ij})\partial _{i}v_{\alpha}\partial _{j}v_{\alpha}\; dv_{\xi}\,.
\]

\n Dividing both sides by $K(n,2)^2\int
_{B(0,\delta)}v_{\alpha}^2\; dv_{\xi}$ and letting $\alpha
\rightarrow +\infty$, we find

\[
(B_0(2,g)+\varepsilon _0)K(n,2)^{-2}\leq K(n,2)^{-2}\limsup \limits
_{\alpha \rightarrow +\infty}\frac{A_{\alpha}}{\int
_{B(0,\delta)}v_{\alpha}^2\; dv_{\xi}}
\]

\begin{equation}\label{new2}
+ \frac{1}{2}\limsup \limits _{\alpha \rightarrow
+ \infty}\frac{B_{\alpha}}{\int _{B(0,\delta)}v_{\alpha}^2\; dv_{\xi}} + \limsup \limits _{\alpha \rightarrow +\infty}\frac{C_{\alpha}}{\int
_{B(0,\delta)}v_{\alpha}^2\; dv_{\xi}},
\end{equation}

\n where

\[
A_{\alpha}=\int _{B(0,\delta)}\eta ^2v_{\alpha}^{2^*}\; dv_{\xi}-\left
( \int _{B(0,\delta)} (\eta v_{\alpha} )^{2^*}\; dv_{\xi}
\right )^{\frac{2}{2^*}},
\]

\[
B_{\alpha}=\int _{B(0,\delta)} (
\partial _k (h_{\alpha}^{ij}\Gamma (h_{\alpha})^k_{ij}
 )+\partial _{ij}h_{\alpha}^{ij} )
(\eta v_{\alpha})^2\; dv_{\xi}
\]

\n and

\[
C_{\alpha}=\left |  \int _{B(0,\delta)}\eta
^2(h_{\alpha}^{ij}-\delta _{ij})\partial _{i}v_{\alpha}\partial
_{j}v_{\alpha}\; dv_{\xi} \right |\ .
\]

\n A simple computation, using the convergence $g_{\alpha} \rightarrow g$
in the $C^2$-topology, gives

\[
\lim \limits _{\alpha \rightarrow +\infty} (
\partial _k (h_{\alpha}^{ij}\Gamma (h_{\alpha})^k_{ij}
 )+\partial _{ij}h_{\alpha}^{ij} )
(0)=\frac{1}{3}Scal_g(x_0),
\]

\n so that, with the step $5$,

\begin{equation}\label{new3}
\limsup \limits _{\alpha \rightarrow +\infty}\frac{B_{\alpha}}{\int
_{B(0,\delta)}v_{\alpha}^2\; dv_{\xi}} = \frac{1}{3}Scal_g(x_0) + \varepsilon
_{\delta},
\end{equation}

\n where $\varepsilon_{\delta} \rightarrow 0$ as $\delta \rightarrow 0$.
\n Using again the convergence in the $C^2$-topology together with some computations, as done in \cite{DjDr}, one finds

\begin{equation}\label{new4}
\limsup \limits_{\alpha \rightarrow + \infty} \frac{A_{\alpha}}{\int_{B(0,\delta)} v_{\alpha}^2\; dv_{\xi}} \leq
\frac{n-4}{12(n-1)} K(n,2)^2 Scal_g(x_0) + \varepsilon_{\delta}
\end{equation}

\n and

\begin{equation}\label{new5}
\limsup \limits_{\alpha \rightarrow + \infty} \frac{C_{\alpha}}{\int
_{B(0,\delta)}v_{\alpha}^2 \; dv_{\xi}} \leq \varepsilon_{\delta}\ .
\end{equation}

\n Putting (\ref{new3}), (\ref{new4}) and (\ref{new5}) into
(\ref{new2}), we obtain, for any $\delta > 0$ small enough,

\begin{equation}\label{new6}
(B_0(2,g) + \varepsilon _0) K(n,2)^{-2} \leq \frac{n-2}{4(n-1)}Scal_g(x_0) + \varepsilon_{\delta}\ .
\end{equation}

\n Letting $\varepsilon_{\delta} \rightarrow 0$ as $\delta \rightarrow 0$
in (\ref{new6}), we arrive at the desired contradiction, since for $n \geq
4$ we have

\[
(B_0(2,g) + \varepsilon_0) K(n,2)^{-2} > \frac{n-2}{4(n-1)} Scal_g(x_0)\ .
\]

\n The $C^2$-topology is sharp as shows the following counter-example. Let $(M,g)$ be a smooth compact Riemannian manifold of dimension $n\geq 4$. Consider a sequence $(f_\alpha)_\alpha \subset C^{\infty }(M)$ of positive functions converging to the constant function $f_0 = 1$ in $L^p(M)$, $p>n$, such that $\max_{M}f_\alpha \rightarrow + \infty$. Let $u_\alpha\in C^{\infty }(M)$, $u_\alpha>0$, be the unique solution of

\[
- \frac{4(n-1)}{n-2} \Delta _{g}u + u = f_\alpha \ .
\]

\n From the classical elliptic $L^p$ theory, it follows that
 $(u_\alpha)_\alpha$ is bounded in $H^{p}_2(M)$, where $H^{p}_2(M)$ stands for the
 second order $L^p$-Sobolev space on $M$, so that $u_\alpha$ converges
 to $u_0$ in $C^{1,\beta}(M)$ for some $0 < \beta < 1$. Moreover,
 $u_0=1$, since $f_\alpha$ converges to $1$ in $L^p(M)$ and the constant
 function $1$ is the unique solution of the limit problem. Therefore,
 $g_\alpha = u_\alpha^{2^*-2}g$ is a smooth Riemannian metric converging to $g$
 in the $C^{1,\beta}$-topology. Note also that there exists a constant
 $c > 0$, independent of $\alpha$, such that

\[
Scal_{g_\alpha} = (- \frac{4(n-1)}{n-2} \Delta _{g}u_\alpha + Scal_{g}
 u_\alpha
 ) u_\alpha^{1-2^*} \geq f_\alpha u_\alpha^{1-2^*} - cu_\alpha^{2-2^*},
\]

\n so that $\max_{M}Scal_{g_\alpha} \rightarrow + \infty$, where
 $Scal_{g}$ denotes the scalar curvature of the metric
$g$. On the other hand, for $n \geq 4$, we have the well-known lower
 bound (see \cite{H3})

\[
B_0(2,g_\alpha) \geq \frac{n-2}{4(n-1)} K(n,2)^2
 \max_{M}Scal_{g_\alpha},
\]

\n so that $B_0(2,g_\alpha)\rightarrow + \infty$. In particular,
 $B_0(2,g_\alpha) \not \rightarrow B_0(2,g)$.

\bl\\

\section{Proof of Theorem \ref{cont1}}

\n Part of the proof in the case $1 < p < 2$ follows a similar outline
 to the proof of Theorem \ref{cont2}. The proof in the case $p = 1$ is inspired in \cite{D4}. So, we will present a resumed
 proof, emphasizing the essential points. Consider initially a sequence
 $(g_\alpha)_{\alpha} \subset {\cal M}_2$ converging to a smooth metric $g$ in the
 $C^0$-topology. Similarly to the previous proof, the convergence in the $C^2$-topology will be used only
 in the last step. As before, we suppose, by contradiction, that the continuity
 fails, so that there exists $\varepsilon_0 > 0$ such that

\[
|B_0(p,g_{\alpha}) - B_0(p,g)| > \varepsilon_0
\]

\n for infinitely many $\alpha$. Then, at least, one of the situations
 holds:

\[
B_0(p,g) - B_0(p,g_{\alpha}) > \varepsilon_0\ \ {\rm or}\ \
 B_0(p,g_{\alpha}) - B_0(p,g) > \varepsilon_0
\]

\n for infinitely many $\alpha$. If the first alternative holds, then
 for any $u \in H^p_1(M)$,

\[
\left( \int_{M} |u|^{p^*}\; dv_{g_{\alpha}} \right)^{p/p^*} \leq
 K(n,p)^p \int_{M} |\nabla_{g_{\alpha}} u|^p \; dv_{g_{\alpha}} +
 (B_0(p,g) -
 \varepsilon_0) \int_{M} |u|^p\; dv_{g_{\alpha}},
\]

\n so that, letting $\alpha \rightarrow +\infty$, we obtain the first
 contradiction. Suppose then that $B_0(p,g) + \varepsilon_0 <
 B_0(p,g_{\alpha})$ for infinitely many $\alpha$. For each $\alpha > 0$, we
 consider
 the functional

\[
J_{\alpha,p}(u)= \int_{M} |\nabla_{g_{\alpha }} u|^p \; dv_{g_{\alpha}}
 +(B_0(p,g) + \varepsilon_0) K(n,p)^{-p} \int_{M}  |u|^p \;
 dv_{g_{\alpha}}
\]

\n defined on $\Lambda_{\alpha,p} =\{ u\in H^p_1(M):\;  \int_{M}
 |u|^{p^*} \; dv_{g_{\alpha}}=1\}$. From the definition of
 $B_0(p,g_{\alpha})$, it follows that

\[
\lambda_{\alpha,p} := \inf_{u \in \Lambda_{\alpha,p}} J_{\alpha,p}(u) <
 K(n,p)^{-p}\ .
\]

\n Assume first $1 < p < 2$. In this case, the inequality above combined with the classical elliptic theory of quasi-linear
 elliptic operators imply the existence of a positive minimizer
 $u_{\alpha}
 \in \Lambda_{\alpha,p}$ for $\lambda_{\alpha,p}$ with $u_\alpha \in
 C^1(M)$. The Euler-Lagrange equation for $u_{\alpha}$ is then

\begin{gather} \label{EL-p}
-\Delta_{p, g_{\alpha}}u_{\alpha} + (B_0(p,g) + \varepsilon_0)
 K(n,p)^{-p}u_{\alpha}^{p-1} = \lambda_{\alpha,p} u_{\alpha }^{p^*-1}\
 .
 \tag{$E_{\alpha,p}$}
\end{gather}

\n Here, $\Delta_{p,g_{\alpha}}u = {\rm
 div}_{g_{\alpha}}(|\nabla_{g_{\alpha}} u|^{p-2} \nabla_{g_{\alpha}} u)$
 denotes the $p$-Laplacian operator
 associated to the metric $g_{\alpha}$. By the $C^0$-convergence of
 $g_\alpha$, there exist constants $\gamma_1, \gamma_2 > 0$ such that
 $\gamma_1 g \leq g_{\alpha} \leq \gamma_2 g$ for $\alpha > 0$ large. This
 fact
 combined with $J_{\alpha,p}(u_{\alpha}) < K(n,p)^{-p}$ imply that
 $(u_{\alpha})_{\alpha}$ is bounded in $H_1^p(M)$, so that $u_{\alpha}
 \rightharpoonup u$ in $H_1^p(M)$, $u \geq 0$, and

\[
\int_{M} u_{\alpha}^q \; dv_{g_{\alpha}} \rightarrow  \int_{M} u^q \;
 dv_{g},
\]

\n for any $1 \leq q < p^*$, as $\alpha \rightarrow +\infty$, up to a
 subsequence. Let $\lambda_{\alpha,p} \rightarrow \lambda_p$. Since $g_\alpha \rightarrow g$ in $C^0$, letting $\alpha
 \rightarrow +\infty$ in (\ref{EL-p}) and using measure theory standard
 arguments, we find that $u$ satisfies in the weak sense

\begin{gather} \label{ELLIM-p}
- \Delta_{p,g} u + (B_0(p,g) + \varepsilon_0) K(n,p)^{-p}u^{p-1} =
 \lambda_p u^{p^*-1}\ . \tag{$E_p$}
\end{gather}

\n If $u \not \equiv 0$, from (\ref{B-opt}),
 (\ref{ELLIM-p}) and $0 \leq \lambda_p \leq K(n,p)^{-p}$, one obtains
 $\int_{M} u^{p^*}\; dv_g > 1$, and this contradicts

\[
\int _M u^{p^*}dv_{g} \leq \liminf_{\alpha \rightarrow + \infty} \int
 _M u_{\alpha }^{p^*}\; dv_{g_{\alpha}} = 1\ .
\]

\n Assume then that $u \equiv 0$ on $M$. Arguing of similar manner to
 the case $p=2$, one gets $\lambda_{\alpha,p} \rightarrow K(n,p)^{-p}$.
 The proofs of the steps from $1$ to $5$, in the $p = 2$ case, relied strongly on local
 H\"{o}lder estimates of weak solutions of elliptic equations and on De
 Giorgi-Nash-Moser type iterative schemes. These tools are also valid
 in
 the quasi-linear elliptic context, we refer to \cite{LU}
 and \cite{Se} for results in the quasi-linear elliptic
 theory. For equations as above, involving a family of $p$-Laplacian
 divergence type operators associated to $g_\alpha$, $\alpha > 0$, such
 tools
 require only $C^0$-convergence of $g_\alpha$. In fact, one needs only
 constants $\gamma, c_0 > 0$, independent of $\alpha$,
 such that $g_\alpha \geq \gamma \xi$ and $||(g_\alpha)_{ij}||_{C^0}
 \leq c_0$ for $\alpha > 0$ large. So, the steps from $1$ to $5$ extend
 readily to $1 < p < 2$, we refer to \cite{D2} for more details.
 Therefore,
 for $1 < p < 2$, these steps take the following form:

We say that $x \in M$ is a point of concentration of $(u_{\alpha})_{\alpha}$ if,
 for any $\delta >0$,

\[
\limsup \limits_{\alpha \rightarrow +\infty} \int_{B_g (x, \delta)}
 u_{\alpha }^{p^*} \; dv_{g_{\alpha}} > 0 \ .
\]

\n {\bf Step 1:} The sequence $(u_{\alpha})_{\alpha}$ possesses exactly one
 point of concentration $x_0$, up to a subsequence.\\

\n {\bf Step 2:} Let $x_0 \in M$ be the unique point of concentration of
 $(u_{\alpha})_{\alpha}$. Then,

\begin{equation}\label{e8}
\lim \limits _{\alpha \rightarrow +\infty} u_{\alpha} = 0\ \ {\rm in}\
 \ C_{loc}^0(M \setminus \{x_0\})\ .
\end{equation}

\n {\bf Step 3:} For each $R > 0$, one has

\begin{equation} \label{first}
\lim _{\alpha \rightarrow  + \infty} \int _{B_{g_{\alpha}}(x_{\alpha},
 R \mu_{\alpha})} u_{\alpha }^{p^*}\; dv_{g_{\alpha }} = 1 -
 \varepsilon_{R}
\end{equation}

\n where $\mu_\alpha = ||u_\alpha||_{\infty}^{- p^*/n}$ and
 $\varepsilon = \varepsilon_R \rightarrow 0$ as $R \rightarrow +
 \infty$.\\

\n {\bf Step 4:} There exists a constant $c > 0$, independent of
 $\alpha$, such that

\[
d_{g}(x,x_{\alpha})^{n/p^*} u_{\alpha}(x) \leq c
\]

\n for all $x \in M$ and $\alpha > 0$ large.\\

\n {\bf Step 5:} For each $\delta > 0$ small, one has

\begin{equation} \label{concent}
\lim_{\alpha \rightarrow +\infty} \frac{\int_{M \setminus
 B_{g}(x_0,\delta)} u_{\alpha}^{p}\; dv_{g}}{\int_{M} u_{\alpha}^p \;
 dv_{g}} = 0\ .
\end{equation}

\n {\bf Step 6:} Here, we assume that $g_\alpha$ converges to $g$ in
the $C^2$-topology. This convergence implies that

\[
\liminf \limits _{\alpha \rightarrow +\infty}
{\rm inj}_{g_{\alpha}}(M)>0\ .
\]

\n Thus, there exists $\delta > 0$ small enough and independent of $\alpha$ such that
$B_{g_{\alpha}}(x_{\alpha},\delta)$ is a geodesic ball for all
$\alpha > 0$ large. In addition, $\exp_{x_{\alpha}, g_\alpha} \circ \exp_{x_0, g}^{-1}$ converges to the identity map $id: B(0, \delta) \rightarrow \R^n$ in the $C^3$-topology. Consider a smooth function $\eta$ such that
$0 \leq \eta \leq 1$, $\eta = 1$ in $(0, \delta)$, $\eta = 0$ in $(2
\delta, + \infty)$ and $|\nabla \eta| \leq c/\delta$ for some
constant $c > 0$ independent of $\delta$. Define $\eta_{\alpha}(x) =
\eta(d_g(x, x_\alpha))$. In what follows, several possibly different
positive constants, independent of $\delta$ and $\alpha$, will be
denoted by $c$. Since $x_{\alpha} \rightarrow x_0$ and $g_\alpha
\rightarrow g$ in the $C^2$-topology, the Cartan expansion
of $g_\alpha$ in a normal coordinates system gives for $\alpha > 0$
large,

\begin{equation} \label{2a}
( 1 - c d_{g_\alpha}(x,x_{\alpha})^2 ) dv_{g_\alpha} \leq dv_{\xi} \leq
 ( 1 + c d_{g_\alpha}(x,x_{\alpha})^2 ) dv_{g_\alpha}\ .
\end{equation}

\n and

\begin{equation} \label{1a}
|\nabla (\eta_{\alpha} u_{\alpha})|^{p}(x) \leq
 |\nabla_{g_\alpha} (\eta_{\alpha} u_{\alpha})|^{p}(x) ( 1 + c
 d_{g_\alpha}(x,x_{\alpha})^2 )
\end{equation}

\n Clearly, (\ref{2a}) gives

\begin{equation} \label{3a}
\int_{B_g(x_{\alpha}, 2\delta)} ( \eta_{\alpha} u_{\alpha} )^{p^*} \;
 dv_{\xi} \geq 1 - \int_{M \setminus B_g(x_{\alpha}, \delta)}
 u_{\alpha}^{p^*} \; dv_{g_\alpha} - c \int_{B_g(x_{\alpha}, 2\delta)}
 u_{\alpha}^{p^*} d_{g_\alpha}(x,x_{\alpha})^2 \; dv_{g_\alpha}\ .
\end{equation}

\n By the step $2$,

\[
\int_{M \setminus B_g(x_{\alpha}, \delta)} u_{\alpha}^{p^*}\;
 dv_{g_\alpha} = o(||u_{\alpha}||_{p}^{p}),
\]

\n and, by the step $4$,

\[
\int_{B_g(x_{\alpha}, 2\delta)} u_{\alpha}^{p^*}
 d_{g_\alpha}(x,x_{\alpha})^2 \; dv_{g_\alpha} \leq c \delta^{2-p}
 ||u_{\alpha}||_{p}^{p}\ .
\]

\n So, (\ref{3a}) yields

\begin{equation}\label{4a}
\left( \int_{B_g(x_{\alpha}, 2\delta)} ( \eta_{\alpha}u_{\alpha}
 )^{p^*}\; dv_{\xi} \right)^{p/p^*} \geq 1 -
 o(||u_{\alpha}||_{p}^{p}) -
 c \delta^{2-p} ||u_\alpha||_{p}^{p}\ .
\end{equation}

\n By (\ref{2a}) and (\ref{1a}), we also have

\begin{equation}\label{5a}
K(n,p)^{p} \int_{B_g(x_{\alpha}, 2\delta)} |\nabla ( \eta_{\alpha}
 u_{\alpha} ) |^{p} \; dv_{\xi} \leq K(n,p)^{p} \int_{B_g(x_{\alpha},
 2\delta)} |\nabla_{g_\alpha} ( \eta_{\alpha} u_{\alpha} )|^{p} \;
 dv_{g_\alpha}
\end{equation}

\[
+ c \int_{B_g(x_{\alpha}, 2\delta)} |\nabla_{g_\alpha} ( \eta_{\alpha}
 u_{\alpha} )|^{p} d_{g_\alpha}(x,x_{\alpha})^2 \; dv_{g_\alpha}\ .
\]

\n Independently, using that $J_{\alpha,p}(u_{\alpha}) =
 \lambda_{\alpha,p}$, $u_{\alpha} \in \Lambda_{\alpha,p}$ and
 $\lambda_{\alpha,p} <
 K(n,p)^{-p}$, one obtains

\[
K(n,p)^{p} \int_{B_g(x_{\alpha}, 2\delta)} |\nabla_{g_\alpha} (
 \eta_{\alpha} u_{\alpha} )|^{p}\; dv_{g_\alpha} \leq 1 - \int_{M
 \setminus
 B_g(x_{\alpha}, \delta)} u_{\alpha}^{p^*} \; dv_{g_\alpha} - (B_0(p,g)
 +
 \varepsilon_0) \int_{B_g(x_{\alpha}, \delta)} u_{\alpha}^{p} \;
 dv_{g_\alpha}
\]

\[
+ c \delta^{-p} \int_{M \setminus B_g(x_{\alpha}, \delta)}
 u_{\alpha}^{p} \; dv_{g_\alpha} + c \int_{M \setminus B_g(x_{\alpha},
 \delta)}
 |\nabla_{g_\alpha} u_{\alpha}|^{p} \; dv_{g_\alpha}\ .
\]

\n In order to estimate the remaining integrals, consider a smooth
 function $\zeta$ such that $0 \leq \zeta \leq 1$, $\zeta = 0$ in $(0,
 \delta)$, $\zeta = 1$ in $(\delta, + \infty)$ and define
 $\zeta_{\alpha}(x)
 = \zeta(d_g(x, x_\alpha))$. Taking $\zeta_{\alpha}^p u_\alpha$ as a
 test function in (\ref{EL}), integrating by parts, using Young's
 inequality, one finds

\[
\int_{M} \zeta_\alpha^{p} |\nabla_{g_\alpha} u_{\alpha}|^{p} \;
 dv_{g_\alpha} \leq \int_{M} \zeta_\alpha^{p} u_{\alpha}^{p^*} \;
 dv_{g_\alpha}
 + c \int_{M \setminus B_g(x_{\alpha}, \delta/2)} u_{\alpha}^{p} \;
 dv_{g_\alpha}
\]

\[
\leq \int_{M \setminus B_g(x_{\alpha}, \delta/2)} u_{\alpha}^{p^*} \;
 dv_{g_\alpha} + c \int_{M \setminus B_g(x_{\alpha}, \delta/2)}
 u_{\alpha}^{p} \; dv_{g_\alpha} = o(||u_{\alpha}||_{p}^{p}),
\]

\n so that

\begin{equation} \label{17}
\int_{M\setminus B_g(x_{\alpha}, \delta)} |\nabla_{g_\alpha}
 u_{\alpha}|^{p}\; dv_{g_\alpha} = o(||u_{\alpha}||_{p}^{p})\ .
\end{equation}

\n Thus,

\begin{equation}\label{6a}
K(n,p)^{p} \int_{B_g(x_{\alpha}, 2\delta)} |\nabla_{g_\alpha}
 (\eta_{\alpha} u_{\alpha} )|^{p} \; dv_{g_\alpha} \leq 1 - (B_0(p,g) +
 \varepsilon_0) \int_{B_g(x_{\alpha}, \delta)} u_{\alpha}^{p} \;
 dv_{g_\alpha} +
 o(||u_{\alpha}||_{p}^{p})
\end{equation}

\n Taking now $\eta_{\alpha}^p u_\alpha d_{g_\alpha}(x,x_{\alpha})^2$
 as a test function in (\ref{EL}), where $\eta_{\alpha}$ is given in
 the
 beginning of this step, integrating by parts and again using Young's
 inequality, one obtains

\[
\int_{B_g(x_{\alpha}, 2\delta)} \eta_\alpha^{p} |\nabla_{g_\alpha}
 u_{\alpha}|^{p} d_{g_\alpha}(x,x_{\alpha})^2 \; dv_{g_\alpha} \leq
 \int_{B_g(x_{\alpha}, 2\delta)} u_{\alpha}^{p^*}
 d_{g_\alpha}(x,x_{\alpha})^2
 \; dv_{g_\alpha} + c \int_{M \setminus B_g(x_{\alpha}, \delta)}
 u_{\alpha} |\nabla_{g_\alpha} u_{\alpha}|^{p - 1} \; dv_{g_\alpha}
\]

\[
+ c \int_{B_g(x_{\alpha}, 2\delta)} d_{g_\alpha}(x,x_{\alpha})^{2 - p}
 \eta_\alpha u_{\alpha} \eta_\alpha^{p - 1} |\nabla_{g_\alpha}
 u_{\alpha}|^{p - 1} d_{g_\alpha}(x,x_{\alpha})^{p-1} \; dv_{g_\alpha}
\]

\[
\leq \int_{B_g(x_{\alpha}, 2\delta)} u_{\alpha}^{p^*}
 d_{g_\alpha}(x,x_{\alpha})^2 \; dv_{g_\alpha} + c \int_{M \setminus
 B_g(x_{\alpha},
 \delta)} u_{\alpha}^{p} \; dv_{g_\alpha} + c \int_{M \setminus
 B_g(x_{\alpha}, \delta)} |\nabla_{g_\alpha} u_{\alpha}|^{p} \;
 dv_{g_\alpha}
\]

\[
+ \frac{1}{2} \int_{B_g(x_{\alpha}, 2\delta)} \eta_\alpha^{p}
 |\nabla_{g_\alpha} u_{\alpha}|^{p} d_{g_\alpha}(x,x_{\alpha})^2 \;
 dv_{g_\alpha}
 + c \int_{B_g(x_{\alpha}, 2\delta)} d_{g_\alpha}(x,x_{\alpha})^{2 - p}
 u_{\alpha}^{p} \; dv_{g_\alpha},
\]

\n so that

\begin{equation} \label{18}
\int_{B_g(x_{\alpha}, 2\delta)}|\nabla_{g_\alpha} (\eta_{\alpha}
 u_{\alpha} )|^{p} d_{g_\alpha}(x,x_{\alpha})^2\; dv_{g_\alpha} \leq c
 \delta^{2-p} ||u_{\alpha}||_{p}^{p} + o(||u_{\alpha}||_{p}^{p}) \ .
\end{equation}

\n Joining (\ref{4a})-(\ref{18}) and the Euclidean Sobolev inequality

\[
\left( \int_{B_g(x_{\alpha}, 2\delta)} ( \eta_{\alpha} u_{\alpha}
 )^{p^*} \; dv_{\xi} \right)^{p/p^*} \leq K(n,p)^{p}
 \int_{B_g(x_{\alpha}, 2\delta)} |\nabla (\eta_{\alpha} u_{\alpha}
 )|^{p} \; dv_{\xi},
\]

\n we obtain

\[
(B_0(p,g) + \varepsilon_0) \int_{B_g(x_{\alpha}, \delta)}
 u_{\alpha}^{p} \; dv_{g_\alpha} \leq c \delta^{2-p}
 ||u_{\alpha}||_{p}^{p} +
 o(||u_{\alpha}||_{p}^{p})\ .
\]

\n Dividing both sides of this inequality by $||u_{\alpha}||_{p}^{p}$,
 letting $\alpha \rightarrow +\infty$ and again using the step 5, one
 finds

\[
(B_0(p,g) + \varepsilon_0) \leq c \delta^{2-p}
\]

\n for all $\delta > 0$ small. This is clearly a contradiction and prove the continuity of $B_0(p,g)$ on $g$ for $1 < p < 2$. If $p = 1$, we fix a sequence $(p_\alpha)_{\alpha} \subset (1, 2)$, $p_\alpha \rightarrow 1$, and, for each $\alpha > 0$, take a minimizer $u_\alpha \in \Lambda_{\alpha, p_\alpha}$ of $J_{\alpha, p_\alpha}$. Since $g_\alpha$ converges to $g$ in the $C^2$-topology, using some ideas of this proof and a key result due to Druet in \cite{D4}, it follows that $(u_\alpha)_{\alpha}$ converges uniformly as $\alpha \rightarrow + \infty$. But this leads directly to a contradiction, since the limit of $(u_\alpha)_{\alpha}$ is an extremal function associated to $B_0(1,g) + \varepsilon_0$. \bl\\

\n {\bf Acknowledgments:}  The first author was partially
supported by Fapemig.\\

 \end{document}